\documentclass[a4paper, 11pt, twoside, final]{scrartcl}

\usepackage{PREAMBLE}
\usepackage{COMMAND}
\usepackage{HYPHENATION}

\begin{document}
\thispagestyle{empty}
\begin{center}
    \rule{\linewidth}{1pt}\\[0.4cm]
    {\sffamily \bfseries \large Ambiguities in one-dimensional phase retrieval from\\
        magnitudes of a linear canonical transform}\\[10pt]
    {\sffamily\footnotesize Robert Beinert}\\[3pt]
    {\sffamily\footnotesize Institut für Numerische und Angewandte Mathematik}\\[-3pt]
    {\sffamily\footnotesize Georg-August-Universität Göttingen}\\
    \rule{\linewidth}{1pt}
\end{center}

\vspace*{10pt}

{\small
    \noindent
    {\sffamily\bfseries Abstract:} Phase retrieval problems occur in a wide range of
    applications in physics and engineering.  Usually, these problems consist in the
    recovery of an unknown signal from the magnitudes of its Fourier transform.  In some
    applications, however, the given intensity arises from a different transformation such
    as the Fresnel or fractional Fourier transform.  More generally, we here consider the
    phase retrieval of an unknown signal from the magnitudes of an arbitrary linear
    canonical transform.  Using the close relation between the Fourier and the linear
    canonical transform, we investigate the arising ambiguities of these phase retrieval
    problems and transfer the well-known characterizations of the solution sets from the
    classical Fourier phase retrieval problem to the new setting.

    \smallskip

    \noindent
    {\sffamily\bfseries Key words:} Phase retrieval; One-dimensional signals; Compact
    support; Linear canonical transform

    \smallskip
    
    \noindent
    {\sffamily\bfseries AMS Subject classifications:} 42A05, 94A08, 94A12
}

\section{Introduction}

In many applications in physics and engineering such as crystallography
\cite{Mil90,Hau91}, astronomy \cite{BS79,DF87}, and laser optics \cite{SST04,SSD+06}, one
is faced with the so-called phase retrieval problem.  The one-dimensional varient of this
problem consists in the recovery of an unknown signal $f \colon \R \rightarrow \C$ in
$L^2(\R)$ with compact support from its Fourier intensity
\begin{equation*}
    \abs{\Fourier \mleft[ f \mright] \mleft( \omega \mright)}
    \coloneqq \absB{ \int_{-\infty}^\infty f \mleft( t \mright) \, \e^{-\I \omega t} \diff
    t}.
    \addmathskip
\end{equation*}
Physically, one can interpret these magnitudes as intensity measurements of a wave in the
far field.  If one measures the intensity in the near field, one has to replace the
Fourier transform by the Fresnel or the fractional Fourier transform \cite{Goo96}.  In
order to investigate the occurring ambiguities of the corresponding phase retrieval
problems, we exploit that all three transformations -- Fourier, Fresnel, and fractional
Fourier transform -- are special cases of the linear canonical transform
\cite[Chap.~9]{Wol79}.

\section{The linear canonical transform}
\label{sec:fresnel-transform}

For the real parameters $a$, $b$, $c$, and $d$ with $ad-bc=1$, the \emph{linear canonical
    transform} is defined by
\begin{equation}
    \label{eq:lin-can-trans}
    \LCT_{(a,b,c,d)} \mleft[ f \mright] \mleft( \omega \mright)
    \coloneqq
    \int_{-\infty}^\infty f \mleft( t \mright) \, K_{(a,b,c,d)} \mleft( \omega, t \mright)
    \diff t
    \submathskip
\end{equation}
with the kernel
\begin{equation}
    \label{eq:kernel}
    K_{(a,b,c,d)} \mleft( \omega, t \mright)
    \coloneqq
    \begin{cases}
        \tfrac{1}{\sqrt{2\pi b}} \, \e^{- \I \frac{\pi}{4}} \,
        \e^{\frac{\I}{2} (\frac{a}{b} t^2 - \frac{2}{b} \omega t + \frac{d}{b} \omega^2)}
        & b \ne 0,
        \\[\fskip]
        \tfrac{1}{\sqrt{a}} \, \e^{ \I \frac{c }{2a} \omega ^2} \,
        \delta \mleft( t - \tfrac{\omega}{a} \mright)
        & \text{else,}
    \end{cases}
    \addmathskip
\end{equation}
where $\delta$ denotes the Dirac delta-function, see \cite[Chap.~9]{Wol79}.

Obviously, the linear canonical transform $\LCT_{(0,1,-1,0)}$ is identical to the Fourier
transform $\Fourier$ up to the multiplicative constant
\begin{equation*}
    \theta \coloneqq \theta_{(a,b,c,d)} \coloneqq \tfrac{1}{\sqrt{2\pi b}} \, \e^{- \I \frac{\pi}{4}}.
\end{equation*}
Moreover, the linear canonical transform covers a complete family of well-known integral
transformations.  For instance, $\LCT_{(1,\nicefrac{1}{2\alpha},0,1)}$ and
$\LCT_{(\cos \alpha, \sin \alpha, - \sin \alpha ,\cos \alpha)}$ with $\alpha \in \R$
coincide with the Fresnel transform \cite{Gor81} and with the fractional Fourier transform
\cite{PD01}, respectively.

If $b \ne 0$, the linear canonical transform can be written in the form
\begin{equation}
    \label{eq:LCT-Fourier}
    \LCT_{(a,b,c,d)} \mleft[ f \mright] \mleft( \omega \mright)
    = \theta_{(a,b,c,d)} \,\e^{ \I \frac{d}{2b} \omega^2}
    \Fourier \mleft[ f \mleft( \cdot \mright) \, \e^{\I \frac{a}{2b} \cdot^2} \mright]
    \mleft( \tfrac{\omega}{b} \mright).
\end{equation}
Using this relation to the Fourier transform, one can easily show that the inverse linear
canonical transform is given by
\begin{equation*}
    \LCT_{(a,b,c,d)}^{-1} [ \lct f ] \mleft( t \mright)
    =
    \int_{-\infty}^\infty \lct f \mleft( \omega \mright) \, \overline{K_{(a,b,c,d)}
        \mleft( \omega, t \mright)} \diff \omega,
\end{equation*}
which coincides up to a unimodular constant with $\LCT_{(d,-b,-c,a)}$.

\section{Phase retrieval from magnitudes of the linear canonical transform}
\label{sec:ambig-fresn-phase-retr}

We now consider the corresponding phase retrieval problem.  In other words, we wish to
recover a signal $f \in L^2(\R)$ with compact support from $\absn{\LCT_{(a,b,c,d)}[f]}$.
Since $\LCT_{(a,0,c,d)}[f]$ is merely a scaled and modulated version of $f$, and since the
recovery of a complex-valued function is not possible from its modulus in general, we
assume that $b \ne 0$.  Similarly to the Fourier setting, this phase retrieval problem
cannot be solved uniquely.

\pagebreak

\begin{Proposition}
    \label{prop:triv-amb-LCT}%
    Let $f \in L^2(\R)$ be a signal with compact support.  Then
    \begin{enumerate}[\upshape(i)]
    \item the rotated signal
        $\e^{\I \alpha} f$ with $\alpha \in \R$
    \item the shifted signal
        $\e^{- \nicefrac{\I a t_0 \cdot}{b}} f(\cdot - t_0)$ with
        $t_0 \in \R$
    \item the reflected signal
        \raisebox{0pt}[0pt][0pt]{$\e^{- \nicefrac{\I a \cdot^2}{b}} \, \overline{f(-\cdot)}$}
    \end{enumerate}
    have the same linear canonical intensity $\absn{\LCT_{(a,b,c,d)} [f]}$.
\end{Proposition}

\begin{Proof}
    The assertion can be established by applying \eqref{eq:LCT-Fourier} and using the
    properties of the Fourier transform.
    \begin{enumerate}[(i)]
    \item 
        $\LCT_{(a,b,c,d)} [ \e^{\I \alpha} f ] ( \omega )
        = \e^{\I \alpha} \LCT_{(a,b,c,d)} [ f ]( \omega )$
    \item
        $\begin{aligned}[t]
            \LCT_{(a,b,c,d)} [ \e^{- \nicefrac{\I a t_0 \cdot}{b}} f(\cdot - t_0) ] (
            \omega )
            &= \theta \, \e^{-\nicefrac{\I a t_0^2}{2b}} \,\e^{\nicefrac{\I d \omega^2}{2b}}
            \Fourier \bigl[ f (\cdot - t_0)\, \e^{\nicefrac{\I a (\cdot-t_0)^2}{2b}} \bigr] (\nicefrac{\omega}{b})
            \\
            &= \e^{-\nicefrac{\I a t_0^2}{2b}} \, \e^{- \nicefrac{\I \omega t_0}{b}}
            \LCT_{(a,b,c,d)} [f] (\omega)
        \end{aligned}$
    \item
        $\begin{aligned}[t]
            \LCT_{(a,b,c,d)} [ \e^{- \nicefrac{\I a \cdot^2}{b}} \overline{f(- \cdot)} ] (
            \omega )
            &= \theta \,\e^{\nicefrac{\I d \omega^2}{2b}}
            \Fourier \bigl[ \overline{f (- \cdot)\, \e^{\nicefrac{\I a (-\cdot)^2}{2b}}}
            \bigr]
            (\nicefrac{\omega}{b})
            \\
            &=  \e^{2\I \arg \theta} \, \e^{\nicefrac{\I d \omega^2}{b}} \,
            \overline{\LCT_{(a,b,c,d)} [f] (\omega)}
        \end{aligned}$
    \end{enumerate}
    Considering the absolute value of each equation finishes the proof.  \qed
\end{Proof}

Without further information about the unknown signal, these three ambiguities cannot be
avoided.  Considering that these signals are, however, closely related to the original
signal $f$, we call them \emph{trivial ambiguities}.  Besides these ambiguities, the phase
retrieval problem usually possesses a series of further non-trivial ambiguities.  Using
\eqref{eq:LCT-Fourier}, the complete solution set can be characterized similarly to the
Fourier case in \cite{Wal63,Hof64,Bei16}.  For this, we denote the Laplace transform and
the autocorrelation function of a signal $f\in L^2 (\R)$ by
\begin{equation*}
    \Laplace \mleft[ f \mright] \mleft( \zeta \mright) \coloneqq \int_{-\infty}^\infty f
    \mleft( t \mright) \, \e^{-\zeta t} \diff t
    \qquad \text{and} \qquad
    A \mleft[ f \mright] \mleft( \zeta \mright)
    \coloneqq \int_{-\infty}^\infty \int_{-\infty}^\infty
    \overline{f \mleft( s \mright)} \, f \mleft( s+t \mright) \, \e^{-\zeta t} \diff s
    \diff t.
    \addmathskip
\end{equation*}

\begin{Theorem}
    \label{the:char-amb-LCT}%
    Let $f \in L^2(\R)$ be a signal with compact support.  Then each signal
    $g \in L^2(\R)$ with compact support and
    $\absn{\LCT_{(a,b,c,d)}[g]} = \absn{\LCT_{(a,b,c,d)}[f]}$ is of the form
    \begin{equation*}
        \Laplace \bigl[ \theta \, \e^{\frac{\I a}{2b} \cdot^2} g\bigr] \mleft( \zeta \mright)
        = C \, \zeta^m \, \e^{\zeta \gamma} \prod_{j=1}^\infty \left( 1 - \tfrac{\zeta}{\eta_j}
        \right) \e^{\frac{\zeta}{\eta_j}},
    \end{equation*}
    where the absolute value $\absn{C}$ and the imaginary part $\Im \gamma$ of the complex
    constants $C$ and $\gamma$ coincide for all signals $g$, and where $\eta_j$ is chosen
    from the zero pair \raisebox{0pt}[0pt][0pt]{$(\xi_j^{\,}, \overline \xi_j)$} of
    \raisebox{0pt}[0pt][0pt]{$A [\theta \, \e^{\nicefrac{\I a \cdot^2}{2b}} f]$}.
\end{Theorem}

\begin{Proof}
    Since
    $\absn{\LCT_{(a,b,c,d)}[f](\omega)} = \absn{ \theta \Fourier [ f \, \e^{\nicefrac{\I a
                \cdot^2}{2b}}] (\nicefrac{\omega}{b})}$,
    we can identify the phase retrieval problem to recover $f$ from
    $\absn{\LCT_{(a,b,c,d)}[f]}$ with the phase retrieval problem to recover
    \raisebox{0pt}[0pt][0pt]{$ \theta \,f \, \e^{\nicefrac{\I a \cdot^2}{2b}}$} from
    \raisebox{0pt}[0pt][0pt]{$\absn{\theta \Fourier [ f \, \e^{\nicefrac{\I a
                    \cdot^2}{2b}}]}$}.
    Hence, the solutions of both problems differ only by
    \raisebox{0pt}[0pt][0pt]{$\theta \,\e^{\nicefrac{\I a \cdot^2}{2b}}$}.  Using
    \cite[Theorem~3.3]{Bei16} to characterize the solutions of the Fourier phase retrieval
    problem, we immediately obtain the assertion.  \qed%
\end{Proof}

\begin{Remark}
    Since we have assumed that the unknown signal $f$ has a compact support, the
    autocorrelation function $A [\theta \, \e^{\nicefrac{\I a \cdot^2}{2b}} f]$ in
    Theorem~\ref{the:char-amb-LCT} is the analytic continuation of
    $\absn{\Fourier [\theta \, \e^{\nicefrac{\I a \cdot^2}{2b}} f]}^2$ from the complex
    axis to the complex plane, \ie\
    \begin{equation*}
        A\mleft[\theta \, \e^{\I \frac{a}{2b} \cdot^2} f\mright] \mleft(\I \omega\mright)
        = \abs{\Fourier \mleft[\theta \, \e^{\I\frac{a}{2b} \cdot^2} f \mright]
            \mleft(\omega\mright)}^2
        \qquad (\omega \in \R),
    \end{equation*}
    see for instance \cite[Proposition~3.2]{Bei16}.  Hence, the required autocorrelation
    function $A [\theta \, \e^{\nicefrac{\I a \cdot^2}{2b}} f]$ is completely encoded in
    the given intensity $\absn{\LCT_{(a,b,c,d)} [f]}$.  \qed
\end{Remark}

\section{Discretization of the linear canonical phase retrieval problem}
\label{sec:discr-line-canon}

To determine a numerical solution, one has to discretize the problem formulation.  For
this purpose, we replace the continuous-time signal $f \colon \R \rightarrow \C$ by a
discrete-time signal $x \colon \Z \rightarrow \C$.  Analogously to the continuous-time
setting, we assume that the signal $x \in \ell^2(\Z)$ has a finite support, which means that
only finitely many signal components $x[n]$ are non-zero.  Discretizing the integral in
\eqref{eq:lin-can-trans}, we define the linear canonical transform of the signal
$x \coloneqq (x[n])_{n\in\Z}$ by
\begin{equation*}
    \LCT_{(a,b,c,d)} \mleft[ x \mright] \mleft( \omega \mright)
    \coloneqq
    \sum_{n\in\Z} x \mleft[ n \mright] \, K_{(a,b,c,d)} \mleft( \omega, n \mright),
\end{equation*}
where $K_{(a,b,c,d)}$ is again the kernel in \eqref{eq:kernel}.

Analogous to \eqref{eq:LCT-Fourier}, the linear canonical transform can be written as
\begin{equation}
    \label{eq:LCT-Fourier:disc}
    \LCT_{(a,b,c,d)} \mleft[ x \mright] \mleft( \omega \mright)
    = \theta_{(a,b,c,d)} \,\e^{ \I \frac{d}{2b} \omega^2}
    \Fourier \mleft[ x \mleft[ \cdot \mright] \, \e^{\I \frac{a}{2b} \cdot^2} \mright]
    \mleft( \tfrac{\omega}{b} \mright)
\end{equation}
whenever $b \ne 0$.  Here $\Fourier$ denotes the discrete-time variant of the Fourier transform
given by
\begin{equation*}
    \Fourier \mleft[ x \mright] \mleft( \omega \mright) \coloneqq
    \sum_{n \in \Z} x \mleft[ n \mright] \, \e^{-\I \omega n}.
    \addmathskip
\end{equation*}
Reversing \eqref{eq:LCT-Fourier:disc}, we notice that the discrete-time linear canonical
transform can be inverted by
\begin{equation*}
    \LCT_{(a,b,c,d)}^{-1} [\lct x] [n] = \int_{-\pi \abs{ \! b \! }}^{\pi \abs{ \! b \!}} \lct x
    \mleft( \omega \mright) \, \overline{K \mleft( \omega, n \mright)} \diff \omega.
    \addmathskip
\end{equation*}

Based on our definitions, the discrete-time variant of the phase retrieval problem can be
stated as follows: recover the unknown discrete-time signal $x \in \ell^2(\Z)$ with finite
support from $\absn{\LCT_{(a,b,c,d)} [x]}$.  For the same reason as before, we assume that
$b\ne 0$.  Adapting the proof of Proposition~\ref{prop:triv-amb-LCT}, we can simply transfer the
three kinds of trivial ambiguities to the discrete-time setting.

\begin{Proposition}
    \label{prop:triv-amb-LCT:disc}%
    Let $x \in \ell^2(\Z)$ be a signal with finite support.  Then
    \begin{enumerate}[\upshape(i)]
    \item the rotated signal
        $\e^{\I \alpha} \, x$ with $\alpha \in \R$
    \item the shifted signal
        $\e^{- \nicefrac{\I a n_0 \cdot}{b}} \, x[\cdot - n_0]$ with
        $n_0 \in \Z$
    \item the reflected signal
        \raisebox{0pt}[0pt][0pt]{$\e^{- \nicefrac{\I a \cdot^2}{b}} \, \overline{x[-\cdot]}$}
    \end{enumerate}
    have the same linear canonical intensity $\absn{\LCT_{(a,b,c,d)} [x]}$.
\end{Proposition}

In order characterize the non-trivial ambiguities, we will exploit the representations of
the non-trivial solutions of the Fourier phase retrieval problem in \cite{BS79,BP15}.
Denoting the support length of the discrete-time signal $x\in\ell^2(\Z)$ by $N$ , we define
the corresponding autocorrelation signal and autocorrelation polynomial by
\begin{equation*}
    a \mleft[ x \mright] \mleft[ n \mright] \coloneqq \sum_{k \in \Z} x \mleft[ k \mright] \, \overline{x
        \mleft[ k+n \mright]}
    \qquad\text{and}\qquad
    P_A \mleft[ x \mright] \mleft( z \mright) \coloneqq z^{N-1} \, \smashoperator{\sum_{n=-N+1}^{N-1}}
    \, a \mleft[ x \mright]\mleft[ n \mright] \, z^n.
\end{equation*}
Since the autocorrelation signal $a[x]$ possesses the support $\{-N+1, \dots, N-1 \}$, the
autocorrelation polynomial $P_A [x]$ is always a well-defined polynomial of degree $2N-2$.

\begin{Theorem}
    \label{the:char-amb-LCT:disc}%
    Let $x \in \ell^2(\Z)$ be a signal with finite support.  Then each signal $y \in \ell^2(\Z)$ with
    finite support and $\absn{\LCT_{(a,b,c,d)}[y]} = \absn{\LCT_{(a,b,c,d)}[x]}$ is of
    the form
    \begin{equation*}
        \Fourier \bigl[ \theta \, \e^{\frac{\I a}{2b} \cdot^2} y\bigr] \mleft( \omega \mright)
        = \e^{\I (\alpha + \omega n_0)}
        \sqrt{\abs{a [\theta \, \e^{ \I\frac{a }{2b} \cdot^2} x] [N-1]}
            \prod_{j=1}^{N-1} \absn{\beta_j}^{-1}}
        \cdot \prod_{j=1}^{N-1} \left( \e^{-\I \omega} - \beta_j \right)
    \end{equation*}
    where $\alpha \in \R$, $n_0 \in \Z$, and $\beta_j$ is chosen from the zero pair
    $(\gamma_j^{\,}, \overline \gamma_j^{\,-1})$ of
    \raisebox{0pt}[0pt][0pt]{$P_A [\theta \, \e^{\nicefrac{\I a \cdot^2}{2b}} x]$} for
    $j=1, \dots, N-1$.
\end{Theorem}

\begin{Proof}
    Similarly to the equivalent continuous-time statement, the relationship
    \eqref{eq:LCT-Fourier:disc} implies that
    $\absn{\LCT_{(a,b,c,d)}[x](\omega)} = \absn{ \theta \Fourier [ x \, \e^{\nicefrac{\I a
                \cdot^2}{2b}}] (\nicefrac{\omega}{b})}$.
    We can thus reduce the considered phase retrieval problem to the recovery of
    \raisebox{0pt}[0pt][0pt]{$ \theta \,x \, \e^{\nicefrac{\I a \cdot^2}{2b}}$} from
    \raisebox{0pt}[0pt][0pt]{$\absn{\theta \Fourier [ x \, \e^{\nicefrac{\I a
                    \cdot^2}{2b}}]}$}.
    Now, the assertion immediately follows from the characterization of the solution set
    of the Fourier phase retrieval problem in \cite[Theorem~2.4]{BP15}. \qed%
\end{Proof}

\begin{Remark}
    Considering the well-known relation 
    \begin{equation}
        \label{eq:auto-LCT:disc}
        \Fourier \bigl[a \bigl[ \theta \, x \, \e^{\I \frac{a}{2b}\cdot^2} \bigr] \bigr]
        = \absb{\Fourier \bigl[ \theta \, x \, \e^{ \I\frac{a}{2b} \cdot^2} \bigr]}^2
        = \abs{C_{(a,b,c,d)}[x]}^2
    \end{equation}
    between the autocorrelation signal and the squared Fourier intensity of a
    discrete-time signal with finite support, see for instance \cite[\p~1173]{BP15}, the
    autocorrelation polynomial
    \raisebox{0pt}[0pt][0pt]{$P_A [\theta \, \e^{\nicefrac{\I a \cdot^2}{2b}} x]$} in
    Theorem~\ref{the:char-amb-LCT:disc} is completely determined by
    $\absn{\LCT_{(a,b,c,d)}[x]}$.

    Moreover, equation \eqref{eq:auto-LCT:disc} shows that $\absn{\LCT_{(a,b,c,d)}[x]}^2$
    as Fourier transform of the autocorrelation signal
    \raisebox{0pt}[0pt][0pt]{$a [\theta \, \e^{\nicefrac{\I a \cdot^2}{2b}} x]$} is a
    non-negative real-valued trigonometric polynomial of degree $N-1$.  Hence, the
    intensity $\absn{\LCT_{(a,b,c,d)}[x]}$ is already uniquely defined by $2N-1$ samples
    at appropriate points in $[-\pi, \pi)$.  \qed
\end{Remark}

\section{Conclusion}
\label{sec:conclusion}

Using the close relation between the linear canonical transform and the Fourier transform,
we have characterized the complete solution set of the phase retrieval problem to recover
a continuous-time or discrete-time signal from the intensity of its linear canonical
transform.  With the same approach, one can transfer most of the uniqueness results for
the classical phase retrieval problem, see for instance
\cite{BFGR76,KST95,RDN13,BP15,Bei16,BP16} and references therein, to the new setting.

\section*{Acknowledgements}

I gratefully acknowledge the funding of this work by the DFG in the framework of the
SFB~755 \bq{Nanoscale photonic imaging} and of the GRK~2088 \bq{Discovering structure in
    complex data: Statistics meets Optimization and Inverse Problems.}

\bibliographystyle{alphadinUK}
{\footnotesize \bibliography{LITERATURE}}

\newcommand{\etalchar}[1]{$^{#1}$}
\begin{thebibliography}{BFGR76}


\providecommand{\url}[1]{\texttt{#1}}
\expandafter\ifx\csname urlstyle\endcsname\relax
  \providecommand{\doi}[1]{doi: #1}\else
  \providecommand{\doi}{doi: \begingroup \urlstyle{rm}\Url}\fi

\bibitem[Bei16]{Bei16}
\textsc{Beinert}, Robert:
\newblock \emph{One-dimensional phase retrieval with additional interference
  measurements}.
\newblock April 2016. --
\newblock Preprint, arXiv:1604.04489v1

\bibitem[BFGR76]{BFGR76}
\textsc{Burge}, R.~E. ; \textsc{Fiddy}, M.~A. ; \textsc{Greenaway}, A.~H.  ;
  \textsc{Ross}, G.:
\newblock The phase problem.
\newblock {In: }\emph{Proceedings of the Royal Society of London. Series A.
  Mathematical Physical \& Engineering Sciences} 350 (1976), pp. 191--212

\bibitem[BP15]{BP15}
\textsc{Beinert}, Robert ; \textsc{Plonka}, Gerlind:
\newblock Ambiguities in one-dimensional discrete phase retrieval from Fourier
  magnitudes.
\newblock {In: }\emph{Journal of Fourier Analysis and Applications} 21 (2015),
  December, No. 6, pp. 1169--1198

\bibitem[BP16]{BP16}
\textsc{Beinert}, Robert ; \textsc{Plonka}, Gerlind:
\newblock \emph{Enforcing uniqueness in one-dimensional phase retrieval by
  additional signal information in time domain}.
\newblock March 2016. --
\newblock Preprint, arXiv:1604.04493v1

\bibitem[BS79]{BS79}
\textsc{Bruck}, Yu.~M. ; \textsc{Sodin}, L.~G.:
\newblock On the ambiguity of the image reconstruction problem.
\newblock {In: }\emph{Optics communications} 30 (1979), September, No. 3, pp.
  304--308

\bibitem[DF87]{DF87}
\textsc{Dainty}, J.~C. ; \textsc{Fienup}, J.~R.:
\newblock Phase retrieval and image reconstruction for astronomy.
\newblock {In: }\textsc{Stark}, Henry (Ed.): \emph{Image Recovery {\upshape:}
  Theory and Application}.
\newblock Orlando (Florida) : Academic Press, 1987, Chapter ~7, pp. 231--275

\bibitem[Goo96]{Goo96}
\textsc{Goodman}, Joseph~W.:
\newblock \emph{Introduction to Fourier Optics}.
\newblock 2nd Edition.
\newblock New York : McGraw-Hill, 1996 (McGraw-Hill Series in Electrical and
  Computer Engineering {\upshape:} Electromagnetics)

\bibitem[Gor81]{Gor81}
\textsc{Gori}, F.:
\newblock Fresnel transform and sampling theorem.
\newblock {In: }\emph{Optics Communications} 39 (1981), November, No. 5, pp.
  293--297

\bibitem[Hau91]{Hau91}
\textsc{Hauptman}, Herbert~A.:
\newblock The phase problem of x-ray crystallography.
\newblock {In: }\emph{Reports on Progress in Physics} 54 (1991), November, No.
  11, pp. 1427--1454

\bibitem[Hof64]{Hof64}
\textsc{Hofstetter}, Edward~M.:
\newblock Construction of time-limited functions with specified autocorrelation
  functions.
\newblock {In: }\emph{IEEE Transaction on Information Theory} 10 (1964), April,
  No. 2, pp. 119--126

\bibitem[KST95]{KST95}
\textsc{Klibanov}, Michael~V. ; \textsc{Sacks}, Paul~E.  ;
  \textsc{Tikhonravov}, Alexander~V.:
\newblock The phase retrieval problem.
\newblock {In: }\emph{Inverse Problems} 11 (1995), No. 1, pp. 1--28

\bibitem[Mil90]{Mil90}
\textsc{Millane}, R.~P.:
\newblock Phase retrieval in crystallography and optics.
\newblock {In: }\emph{Journal of the Optical Society of America A} 7 (1990),
  March, No. 3, pp. 394--411

\bibitem[PD01]{PD01}
\textsc{Pei}, Soo-Chang ; \textsc{Ding}, Jian-Jiun:
\newblock Relations between fractional operations and time-frequency
  distributions, and their applications.
\newblock {In: }\emph{IEEE Transactions on Signal Processing} 49 (2001),
  August, No. 8, pp. 1638--1655

\bibitem[RDN13]{RDN13}
\textsc{Raz}, Oren ; \textsc{Dudovich}, Nirit  ; \textsc{Nadler}, Boaz:
\newblock Vectorial phase retrieval of 1-D signals.
\newblock {In: }\emph{IEEE Transactions on Signal Processing} 61 (2013), April,
  No. 7, pp. 1632--1643

\bibitem[SSD{\etalchar{+}}06]{SSD+06}
\textsc{Seifert}, Birger ; \textsc{Stolz}, Heinrich ; \textsc{Donatelli}, Marco
  ; \textsc{Langemann}, Dirk  ; \textsc{Tasche}, Manfred:
\newblock Multilevel {G}auss-{N}ewton methods for phase retrieval problems.
\newblock {In: }\emph{Journal of Physics. A. Mathematical and General} 39
  (2006), No. 16, pp. 4191--4206

\bibitem[SST04]{SST04}
\textsc{Seifert}, Birger ; \textsc{Stolz}, Heinrich  ; \textsc{Tasche},
  Manfred:
\newblock Nontrivial ambiguities for blind frequency-resolved optical gating
  and the problem of uniqueness.
\newblock {In: }\emph{Journal of the Optical Society of America B} 21 (2004),
  May, No. 5, pp. 1089--1097

\bibitem[Wal63]{Wal63}
\textsc{Walther}, Adriaan:
\newblock The question of phase retrieval in optics.
\newblock {In: }\emph{Optica Acta: International Journal of Optics} 10 (1963),
  No. 1, pp. 41--49

\bibitem[Wol79]{Wol79}
\textsc{Wolf}, Kurt~B.:
\newblock \emph{Integral Transforms in Science and Engineering}.
\newblock New York : Plenum Press, 1979

\end{thebibliography}

\end{document}